Linearity of Data and Linear Probability Space                11 March 2019
        Companion paper 1

Christopher M Rembold MD
Cardiovascular Division
Department of Internal Medicine
University of Virginia Health System
Charlottesville, Virginia 22908-0146 USA

Running Head: Data Linearity

Contact information:
Christopher M. Rembold, M.D.
Box 800146
Cardiovascular Division
University of Virginia Health System
Charlottesville, Virginia 22908-0146 USA
Telephone (434) 924-2825
Email: crembold@virginia.edu

**ABSTRACT**

Some data is linearly additive, other data is not.  In this paper, I discuss types of data based on the boundedness of the data and their linearity.  1) Unbounded data can be linear.  2) One-side bounded data is usually log transformed to be linear.   3) Two-side bounded data is not linear.   4) Untidy data do not fit in these categories.  An example of two-sided bounded data is probabilities which should be transformed into a linear probability space by taking the logarithm of the odds ratio ($\log_{10}$ odds) which is termed Weight (W).  Calculations of means and standard deviation is more accurate when calculated as W values than when calculated as probabilites.  A methods to analyze untidy data is discussed.

# TYPES OF DATA AND THEIR LINEARITY

Data are data, scientists attempt to measure as accurately as possible. However, once we decide to mathematically look at data, we need to determine whether the data are linear prior to analyzing data. Linearity of data has been discussed, see (*1*)-p. 423, however, it has not been applied widely. I propose a classification of data linearity based on how the data are bounded, i.e. 1) unbounded, 2) one-side bounded, 3) two-side bounded, and 4) untidy.

1) Unbounded data in one dimension goes from negative infinity to positive infinity. A mathematical example is real numbers. Newtonian physics is unbounded and linear: the three dimensions x, y and z are unbounded. Mathematically, the Newtonian time dimension t is treated as if it were unbounded (t is only limited by entropy). Calculation in Newtonian space are linear including calculation of force and acceleration.

Rotational data are also unbounded. The θ in $re^{i\theta}$ is clearly unbounded and although r is traditionally always positive, a negative r occurs with $re^{i\theta} = -re^{i(\theta+\pi)}$ so rotational data are linear.

Calculations in Newtonian and rotational space have been shown to be quite accurate at low velocities. In seeking to account for nonlinearities in Maxwell's electromagnetic theory, Einstein disrupted the linearity of Newtonian physics by showing how it was nonlinear at higher velocities (*2*) and then finding a linear space for calculations (*3*) which he termed spacetime.

Integer data that can be negative and are not to be transformed into probabilities can be unbounded and linear.

2) One-side bounded data goes from a single value to typically positive infinity. The typical single value is a limit approaching 0. An example is the frequency (or wavelength) of a sinusoidal wave. Calculations would be linear if done as a logarithm of the frequency (or wavelength).

Integer data that cannot be negative and are not to be transformed into probabilities could be also be one-side bounded. This is more difficult to imagine.

3) Two-side bounded data exists between two defined values. An example is probabilities which go from 0 to 1 (i.e. 0 to 100%). Probabilities of different events are linearly additive. However, changes in probabilities are not linearly additive and need to be transformed into a linear probability space as the $\log_{10}$ Odds for calculations (see below and Rembold, main paper). Examples include a lot of chemical and biological data such as covalent modification and reversible binding. This could also include changes in the probability of spin direction in quantum mechanics.

Integer data that are to be converted into probabilities, like measuring recurrence of cancer or cardiac ischemia after an intervention, need to be transformed into linear probability space as described for probabilities (see Rembold, main paper).

4) Untidy data. An example is data presented as a percent of untreated control. These data could be one-side unbounded. A rational approach is to consider the characteristics of the most likely underlying process and transform data appropriately. This is will be discussed in detail below.

In this paper, I will discuss probability data in more detail. I show that when data are compared, calculating means and standard deviation is more accurate when calculated in linear probability space than when calculated as probabilities. I propose a method to bring untidy data that is not percentages into linear probability space for testing.

**LINEAR PROBABILITY SPACE**

As described above, probability space has data ranging from 0 to 1 (0 to 100%). As such probability space is two-sided unbounded. Changes in probabilities are not additive. In a companion paper (Rembold, main paper), I suggest calculation of statistics is best done in a linear probability space by calculating as the log of odds. The log odds has been described as "weight of evidence (*4, 5*)," "Weight ("W") (*6*), and also as the logit transform (as the natural log, (*1*)-p. 394). For various reasons, the $\log_{10}$ odds is more intutative (Rembold, main paper). In linear probability space there are three statistics:

1) W (weight) contains information of the probability of something. However, unlike probabilities, changes in W values are additive. The conversions between probability, odds and W are:

$$W = \text{Weight} = \log_{10} \text{Odds} = \log_{10} \frac{\text{Probability}}{1-\text{Probability}} = \log_{10} \frac{\text{Percentage}/100}{1-\text{Percentage}/100}$$

In Excel: W = LOG10(Probability / (1–Probability))     (1)

$$\text{Probability} = \frac{10^W}{1+10^W}$$
In Excel: Probability=POWER(10,W)/(1+POWER(10,W)) (2)

2) I (impact) is a quantity that changes a probability. When one object impacts another, it moves the second object. Therefore a test or treatment can "impact" W by a quantitative measure I. I values are the same as Effect Size which quantitative measures of how much means of two measures differ rather than how certain we are about the difference in means. Within linear probability space, I is easily calculated as the difference of the W values (means) measured in two groups or with two treatments:

Impact = I = $W_{mean2} - W_{mean1}$     (3)

In linear probability space Bayes' theorem reduces to:
$W_{posttest} = W_{pretest} + I$     (4)

3) C (certainty) replaces "significance" to describe confidence about a result. C values are calculated with p values the same way as Ws are calculated with probabilities. Note that a more negative C value implies more certainty about a result. E.g. a C value of –4 corresponds to a p value of ~0.0001 while a C value of –2 corresponds to a p value of ~0.01.

$$C = \text{Certainty} = \log_{10} \frac{\text{p value}}{1-\text{p value}} \quad \text{In Excel: W = LOG10(p value / (1–p value))} \quad (5)$$

**CALCULATION OF PROBAILITIES IS NOT ACCURATE WHEN COMPARED TO LINEAR PROBABLITY SPACE**

Fig. 1 demonstrates the quantitative divergence of mean and standard deviation (SD) at varying means when calculations are done as regular probabilities (percent) vs. linear probability space. I created three data sets with means of 50% and SD of 6.8%, 3.4%, and 1.7% (data shown in the bottom of Fig. 1). I then multiplied the individual data by numbers between 0.00005 to 1 to create data with means between 0.0001% and 50%. Then means and SD were calculated with probabilities and in linear probability space (W). Then I calculated the percent difference in the mean and SD when calculated as probability as compared to calculations as W. Means were the same only at 50% and diverged more at lower mean percentages and with higher SD (highest SD is dashed line, Fig. 1, top panel). SD were never the same, and diverged in the opposite direction as means; the divergence was more as lower mean and with higher SD (Fig. 3, bottom panel). The percent deviation of the SD was about tenfold greater than the percent deviation of the mean. These data suggest when data are compared, calculations should be done in a linear probability space.

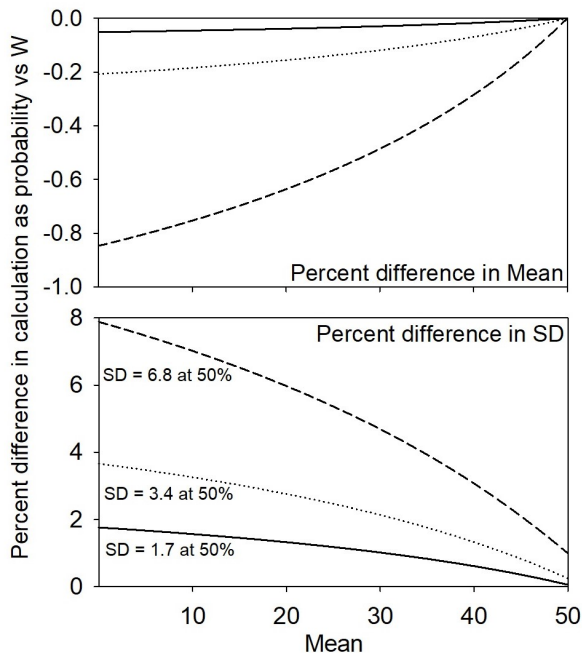

Figure 1. The difference in the mean (top) and SD (bottom) when calculations were done in regular probabilities as compared to calculation done in linear probability space (W) at different mean (x axis in percentages) with 3 different SD (different curves). Data are shown at bottom.

Data for SD 6.8 are 38, 44, 48, 49, 50, 51, 52, 56, 62
Data for SD 3.4 are 44, 47, 49, 49.5, 50, 50.5, 51, 53, 56
Data for SD 1.7 are 47, 48.5, 49.5, 49.75, 50, 50.25, 50.5, 51.5, 53

**AN EXAMPLE WITH "UNTIDY" DATA THAT ARE NOT LINEAR**

Untidy data are data that difficult to linearize. One possibility is to follow the biology or the chemistry or the physics of what is being studied to determine the characteristics of the process and the data and how to linearize it. A biologic example of a probability is the percent phosphorylation of a protein on a single residue. This is stoichiometry with percent phosphorylation ranging from 0% to 100%. These data can be easily converted to linear probability space with Eqns. 1-2 above (see examples in Rembold, main paper). However, sometimes measures of protein phosphorylation are presented as percent of control. Since the underlying process is covalent modification (phosphorylation), then the data should be converted to a linear probability space.

I propose that untidy data needs a transform prior to conversion into linear probability space. Specifically, data can be transformed in three arbitrary ways with 1) no transform, 2) ln x, and 3) $e^x$. Then the mean and SD of all the data can be calculated, then data transformed such that whole data set has mean = 0, mean–1SD = –1, and mean+1SD = +1 (i.e. y = (x – total mean)/total SD). Then probability can be calculated from the SD with Eqn. (6) - this assumes SD is from a normal distribution (this is the major assumption regarding data distribution). Then W can be calculated from W = $\log_{10}$(probability/(1–probability)) (Eqn. 1). The three transforms can be analyzed to see which is more evenly distributed with a skewness function (e.g. (mean – median)/SD). Then the transform with the lowest skewness can be used as the best transform. Each group can then mapped to linear probability space, then can be tested producing I and C values. If needed the result can be inverse mapped back to the original data space.

I propose that these three transforms would not count as a loss of a DF (degree of freedom). If none of the transformations are thought to be adequate, then an investigator chosen function can be substituted with a loss of a DF for each and every constant and operator. If there are values of 0 or 100% which map to W of -∞ and +∞; the solution is to change the 0% and 100% value to 50%/n and 100%–50%/n, respectively, where n is the number of measurements.

The conversions between SD for a normal distribution and probability is

Probability = ½ (1+erf($\frac{SD}{\sqrt{2}}$))     where erf if the error function

in Excel: Probability = ½*(1+ERF(SD/SQRT(2)))                (6)

SD = $\sqrt{2}$*inverf(2·Probability – 1)    where inverf is the inverse error function

in Excel there is no inverse error function so the inverse gamma is used:
SD = IF(Probability>=0.5, SQRT(2)*SQRT(GAMMAINV(2*Probability–1, 0.5, 1)),
   –SQRT(2)*SQRT(GAMMAINV(2*(1–Probability) –1, 0.5, 1)))                (7)

Fig. 2 demonstrates statistical testing of an "untidy" dataset that is NOT easily converted to linear probability space. It shows measurements a protein phosphorylation that can only be calculated as a percent of control, so therefore phosphorylation data are not stoichiometric.

|  | Raw Data | | Normalize x SD | | Prob from SD | | Weight from Prob | | Impact Dif means | 95% CI of Impact | | SD | SE | t SE mult | DF | p (t distr) | C Certainty W of p |
|---|---|---|---|---|---|---|---|---|---|---|---|---|---|---|---|---|---|
|  | 0 min | 10 min | 0 min | 10 min | 0 min | 10 min | 0 min | 10 min |  |  |  |  |  |  |  |  |  |
| Data | 1.00 | 1.29 | -0.72 | -0.46 | 0.24 | 0.32 | -0.51 | -0.32 |  |  |  |  |  |  |  |  |  |
|  | 1.00 | 2.60 | -0.72 | 0.72 | 0.24 | 0.76 | -0.51 | 0.51 |  |  |  |  |  |  |  |  |  |
|  | 1.00 | 2.54 | -0.72 | 0.66 | 0.24 | 0.75 | -0.51 | 0.47 |  |  |  |  |  |  |  |  |  |
|  | 1.00 | 3.98 | -0.72 | 1.95 | 0.24 | 0.97 | -0.51 | 1.58 |  |  |  |  |  |  |  |  |  |
| Mean x | 1.00 | 2.60 | -0.72 | 0.72 | 0.24 | 0.70 | -0.51 | 0.56 | 1.07 | 1.83 | 0.30 | 0.55 | 0.39 | 2.73 | 6 | 0.03397 | -1.45 |
| SD (stdev) | 0.00 | 1.10 | 0.00 | 0.98 | 0.00 | 0.27 | 0.00 | 0.78 |  |  |  |  |  |  |  |  | Marginally Different |
| n | 4 | 4 | 4 | 4 | 4 | 4 | 4 | 4 |  |  |  |  |  |  |  |  |  |
| Total Mean |  | 1.80 |  | 0.00 |  | 0.47 |  | 0.02 |  |  |  |  |  |  |  |  |  |
| Total Std |  | 1.12 |  | 1.00 |  | 0.31 |  | 0.76 |  |  |  |  |  |  |  |  |  |
| Skew |  |  |  |  |  |  |  | 1.40 |  |  |  |  |  |  |  |  |  |
| ln(Data) | 0.00 | 0.25 | -0.78 | -0.34 | 0.22 | 0.37 | -0.56 | -0.23 |  |  |  |  |  |  |  |  |  |
|  | 0.00 | 0.96 | -0.78 | 0.92 | 0.22 | 0.82 | -0.56 | 0.66 |  |  |  |  |  |  |  |  |  |
|  | 0.00 | 0.93 | -0.78 | 0.88 | 0.22 | 0.81 | -0.56 | 0.63 |  |  |  |  |  |  |  |  |  |
|  | 0.00 | 1.38 | -0.78 | 1.68 | 0.22 | 0.95 | -0.56 | 1.31 |  |  |  |  |  |  |  |  |  |
| Mean x | 0.00 | 0.88 | -0.78 | 0.78 | 0.22 | 0.74 | -0.56 | 0.59 | 1.15 | 1.77 | 0.53 | 0.45 | 0.32 | 3.64 | 6 | 0.01086 | -1.96 |
| SD (stdev) | 0.00 | 0.47 | 0.00 | 0.83 | 0.00 | 0.25 | 0.00 | 0.63 |  |  |  |  |  |  |  |  | Marginally Different |
| n | 4 | 4 | 4 | 4 | 4 | 4 | 4 | 4 |  |  |  |  |  |  |  |  |  |
| Total Mean |  | 0.44 |  | 0.00 |  | 0.48 |  | 0.02 |  |  |  |  |  |  |  |  |  |
| Total Std |  | 0.56 |  | 1.00 |  | 0.33 |  | 0.74 |  |  |  |  |  |  |  |  |  |
| Skew |  |  |  |  |  |  |  | 0.88 |  |  |  |  |  |  |  |  |  |
| exp(Data) | 2.72 | 3.62 | -0.52 | -0.47 | 0.30 | 0.32 | -0.36 | -0.33 |  |  |  |  |  |  |  |  |  |
|  | 2.72 | 13.52 | -0.52 | 0.10 | 0.30 | 0.54 | -0.36 | 0.07 |  |  |  |  |  |  |  |  |  |
|  | 2.72 | 12.73 | -0.52 | 0.05 | 0.30 | 0.52 | -0.36 | 0.04 |  |  |  |  |  |  |  |  |  |
|  | 2.72 | 53.62 | -0.52 | 2.39 | 0.30 | 0.99 | -0.36 | 2.07 |  |  |  |  |  |  |  |  |  |
| Mean x | 2.72 | 20.87 | -0.52 | 0.52 | 0.30 | 0.59 | -0.36 | 0.46 | 0.82 | 1.89 | -0.24 | 0.77 | 0.54 | 1.52 | 6 | 0.17932 | -0.66 |
| SD (stdev) | 0.00 | 22.29 | 0.00 | 1.27 | 0.00 | 0.28 | 0.00 | 1.09 |  |  |  |  |  |  |  |  | Indeterminate |
| n | 4 | 4 | 4 | 4 | 4 | 4 | 4 | 4 |  |  |  |  |  |  |  |  |  |
| Total Mean |  | 11.80 |  | 0.00 |  | 0.45 |  | 0.05 |  |  |  |  |  |  |  |  |  |
| Total Std |  | 17.52 |  | 1.00 |  | 0.24 |  | 0.84 |  |  |  |  |  |  |  |  |  |
| Skew |  |  |  |  |  |  |  | 2.57 |  |  |  |  |  |  |  |  |  |

Fig. 2. Statistical testing in a dataset that is NOT easily converted to linear probability space. It shows measurements of the percent of control increase of tyrosine 118 (Y118) paxillin phosphorylation present in swine carotid arterial smooth muscle homogenates frozen prior to and 10 min after depolarization with 109 mM $K^+$ (7). Paxillin was separated on two separate SDS gels and immunostained with antibodies to total paxillin (one gel) and Y118 phosphorylated paxillin (second gel), and the change in phosphorylation is presented as a ratio of immunostained Y118 to total paxillin, so therefore this phosphorylation data are not stoichiometric.

These data were NOT percentages (Fig. 2, column 2 & 3 top section), so the total mean and total SD of all the data were calculated (1.80 and 1.12, respectively). Then the number of SDs from the mean was calculated as y = (x – total mean)/total SD (so now the mean and SD for all data are 0 and 1.0, respectively, Fig. 2, top section, columns 4 & 5). Then probability was calculated from SD with Eqn. (6) - this assumes SD is from a normal distribution (Fig. 2, top section, columns 6 & 7). Then W was calculated from W = ln (probability/(1-probability), top sections, Fig. 2, columns 8 & 9).

Then the raw data was then transformed in two ways: y = ln(x) (center section) and y = exp(x) (bottom section, columns 2 & 3), then SD (columns 4 & 5), probability (columns 6 & 7), and W (columns 8 & 9) were calculated as above. The skewness of the total data for the raw data and each of the transforms was then calculated and the lowest skewness chosen as the appropriate test (shown is Excel's "skewness" function but also could be skew = (mean – median) / SD)).

Statistical testing was then evaluated. In this example skewness was the lowest (0.88) with the ln(x) transform (bolded center section). The dataset show with a C = –1.96 that the mean Y118 paxillin phosphorylation increased from 0 to 10 min of $K^+$ depolarization. This fits into the marginally different certainty category (Rembold, main paper). In probability space, I = 1.15 corresponding to means in the original data of 1 to 2.6 (a very high impact). With either no transform (top) or exp(x) (bottom), the C was closer to zero, i.e. less certain).

A spreadsheet for these calculations is included in the appendix.

This method for transforming non probability data into linear probability space has a number of assumptions and disadvantages, especially calculating probability from SD assuming a normal distribution. I invite statisticians to find a better method.

**DISCUSSION**

There are other options for looking at data. Some, including the physics community, feel that it would be ideal if continuous data presented itself in Standard Deviation (SD) or Standard Error (SE) units. The physics community frequently reports SD multiples (the number of SD separating the means) as a measure of "significance." A more appropriate measure would be SE multiples (the number of SE separating the means). SE multiples is the same as Fischer's t statistic. If data were to present in SD or SE units, statistical difference between two groups can be estimated by the degrees of freedom and the number of SEs separating the means. Unfortunately, 1) each data set has its own SD and SE based on the data. 2) The SD calculated from the data differs from the SD based on a normal distribution. 3) Data sets do not simply presents as multiples of SEs, so a transform is needed. 4) SD and SEs are not additive, so calculations introduce error.

Some data presents itself naturally as probabilities, e.g. as percentages. Probabilities are the most intuitive and some data presents itself naturally as probabilities. Unfortunately, 1) probabilities do not allow simple statistical testing. 2) Calculation of Bayes' theorem is complex. And 3) changes in probabilities are also not linearly additive so calculations introduce error in the mean and SD that increase as means diverge from 50% (Fig. 1). Therefore, transforms to a linear space should be done prior to statistical testing (Fig. 2). Typically these transforms are not done resulting in potentially deceiving results.

Some may suggest that a linear probability space is not appropriate given that the distribution of the entire sample is not known (this is the same rationale by some for not using the t test or Z scores). All issues on data distribution apply equally regardless of whether calculations are done in linear probability space or not.

Overall, I propose that calculation in linear probability space is more accurate and allows easier testing than percentage space.